\documentclass[a4paper,11pt]{amsart}
\usepackage{graphicx}
\usepackage{mathptmx}
\usepackage{amsmath}
\usepackage{amssymb}
\usepackage{enumitem}
\usepackage{xcolor}

\newmuskip\pFqmuskip

\newcommand*\pFq[6][8]{%
  \begingroup 
  \pFqmuskip=#1mu\relax
  \mathcode`=\string"8000
  \begingroup\lccode`\~=`\,
  \lowercase{\endgroup\let~}\pFqcomma
  F^{#2}_{#3}{\left(\genfrac..{0pt}{}{#4}{#5}\bigg|#6\right)}%
  \endgroup
}
\newcommand{\pFqcomma}{\mskip\pFqmuskip}

\newtheorem{theorem}{Theorem}
\newtheorem{lemma}[theorem]{Lemma}
\newtheorem{corollary}[theorem]{Corollary}
\newtheorem{proposition}[theorem]{Proposition}

\begin{document}

\title[Identities on poly-Dedekind Sums]{Identities on poly-Dedekind Sums}

\author{Taekyun  Kim}
\address{Department of Mathematics, Kwangwoon University, Seoul 139-701, Republic of Korea}
\email{tkkim@kw.ac.kr}

\author{Dae San  Kim}
\address{Department of Mathematics, Sogang University, Seoul 121-742, Republic of Korea}
\email{dskim@sogang.ac.kr}

\author{Hyunseok Lee}
\address{ Department of Mathematics, Kwangwoon University, Seoul 139-701, Republic of Korea}
\email{luciasconstant@kw.ac.kr}

\author{Lee-Chae Jang$^{*}$}
\address{Graduate School of Education, Konkuk University, Seoul, Republic of Korea}
\email{lcjang@konkuk.ac.kr}

\thanks{* is corresponding author}

\maketitle

\begin{abstract}
Dedekind sums occur in the transformation behaviour of the logarithm of the Dedekind eta-function under substitutions from the modular group. In 1892, Dedekind showed a reciprocity relation for the Dedekind sums. Apostol generalized Dedekind sums by replacing the first Bernoulli function appearing in them by any Bernoulli functions and derived a reciprocity relation for the generalized Dedekind sums. In this paper, we consider poly-Dedekind sums which are obtained from the Dedekind sums by replacing the first Bernoulli function by any type 2 poly-Bernoulli functions of arbitrary indices and prove a reciprocity relation for the poly-Dedekind sums.
\end{abstract}

\section{\bf Introduction}
In order to give concise definition of the Dedekind sums, we introduce the notation
\begin{equation}
(\!(x)\!)=\left\{\begin{array}{ccc}
    x-[x]-\frac{1}{2}, & \textrm{if $x\notin\mathbb{Z}$,} \\
    0, & \textrm{if $x\in\mathbb{Z}$},
\end{array}\right.  \quad (\mathrm{see}\ [1,3])\label{1}
\end{equation}
 where $[x]$ denotes the greatest integer not exceeding $x$. \par
 It is well known that the Dedekind sums are defined by
 \begin{equation}
 S(h,m)=\sum_{\mu=1}^{m}\bigg(\!\bigg(\frac{\mu}{m}\bigg)\!\bigg)\bigg(\!\bigg(\frac{h\mu}{m}\bigg)\!\bigg),\quad(\mathrm{see}\ [1,3,6-8,11-13], \label{2}
 \end{equation}
 where $h$ is any integer. \par
 From \eqref{2}, we note that
 \begin{equation}
 S(h,m)=\sum_{\mu=1}^{m}\bigg(\frac{\mu}{m}-\frac{1}{2}\bigg)\bigg(\!\bigg(\frac{h\mu}{m}\bigg)\!\bigg)=\sum_{\mu=1}^{m}\frac{\mu}{m}\bigg(\!\bigg(\frac{h\mu}{m}\bigg)\!\bigg),\quad(\mathrm{see}\ [6,7]). \label{3}
 \end{equation}
As is well known, the Bernoulli polynomials are given by
\begin{equation}
\frac{t}{e^{t}-1}e^{xt}=\sum_{n=0}^{\infty}B_{n}(x)\frac{t^{n}}{n!},\quad(\mathrm{see}\ [1-13]).\label{4}
\end{equation}
When $x=0$, $B_{n}=B_{n}(0),\ (n\ge 0)$, are called the Bernoulli numbers. \par
From \eqref{4}, we note that
\begin{equation}
B_{n}(x)=\sum_{l=0}^{n}\binom{n}{l}B_{l}x^{n-l}=(B+x)^{n},\quad (n\ge 0),\quad (\mathrm{see}\ [2-7]),\label{5}
\end{equation}
with the usual convention about replacing $B^{n}$ by $B_{n}$. \par
We observe that
\begin{equation}
\sum_{l=0}^{n-1}e^{lt}=\frac{t}{t(e^{t}-1)}\big(e^{nt}-1\big)=\sum_{j=0}^{\infty}\bigg(\frac{B_{j+1}(n)-B_{j+1}}{j+1}\bigg)\frac{t^{j}}{j!},\quad(n\in\mathbb{N}). \label{6}
\end{equation}
Thus, by \eqref{6}, we get
\begin{equation}
\sum_{l=0}^{n-1}l^{j}=\frac{1}{j+1}\big(B_{j+1}(n)-B_{j+1}\big),\quad (n\in\mathbb{N},\ j\ge 0).\label{7}
\end{equation}
Recently, Kim-Kim considered the polyexponential function of index $k$ given by
 \begin{equation}
 \mathrm{Ei}_{k}(x)=\sum_{n=1}^{\infty}\frac{x^{n}}{n^{k}(n-1)!},\ (k\in\mathbb{Z}),\quad (\mathrm{see}\ [5,10]).\label{8}
 \end{equation}
Note that $\mathrm{Ei}_{1}(x)=e^{x}-1$.\par
In [5], the type 2 poly-Bernoulli polynomials of index $k$ are defined in terms of the polyexponential function of index $k$ as
\begin{equation}
\frac{\mathrm{Ei}_{k}(\log(1+t))}{e^{t}-1}e^{xt}=\sum_{n=0}^{\infty}B_{n}^{(k)}(x)\frac{t^{n}}{n!},\quad(k\in\mathbb{Z}). \label{9}
\end{equation}
When $x=0$, $B_{n}^{(k)}=B_{n}^{(k)}(0)$, $(n\ge 0)$, are called the type 2 poly-Bernoulli numbers of index $k$. Note here that $B_n^{(1)}(x)=B_n(x)$ are the Bernoulli polynomials.\par
The fractional part of $x$ is denoted by
\begin{equation}
\langle x\rangle=x-[x]. \label{10}
\end{equation}
The Bernoulli functions are defined by
\begin{equation}
\overline{B}_{n}(x)=B_{n}(\langle x\rangle),\quad (n\ge 0),\quad(\mathrm{see}\ [1,3,12]). \label{11}
\end{equation}
Thus, by \eqref{3} and \eqref{11}, we get
\begin{align}
S(h,m)\ &=\ \sum_{\mu=1}^{m-1}\frac{\mu}{m}\bigg(\frac{h\mu}{m}-\bigg[\frac{h\mu}{m}\bigg]-\frac{1}{2}\bigg) \label{12} \\
&=\ \sum_{\mu=1}^{m-1}\frac{\mu}{m}\overline{B}_{1}\bigg(\frac{h\mu}{m}\bigg)\ =\ \sum_{\mu=1}^{m-1}\overline{B}_{1}\bigg(\frac{\mu}{m}\bigg)\overline{B}_{1}\bigg(\frac{h\mu}{m}\bigg),    \nonumber
\end{align}
where $h,m$ are relatively prime positive integers. \par

We need the following lemma which is well-known or easily shown.

\begin{lemma}
Let $n$ be a nonnegative integer, and let $d$ be a positive integer. Then we have \\
(a) $\sum_{i=0}^{d-1} B_n\big(\frac{x+i}{d}\big)=d^{1-n}B_n(x)$, \\
(b) $\sum_{i=0}^{d-1} \overline{B}_n\big(\frac{x+i}{d}\big)=d^{1-n}\overline{B}_n(x)$, \\
(c) $\sum_{i=0}^{d-1} B_n\big(\frac{\langle x \rangle+i}{d}\big)=\sum_{i=0}^{d-1} \overline{B}_n\big(\frac{x+i}{d}\big)$, \, \textnormal{for all real} $x$.
\end{lemma}

Dedekind showed that the quantity $ S(h,m)\ =\   \sum_{\mu=1}^{m-1}\frac{\mu}{m}\overline{B}_{1}\big(\frac{h\mu}{m}\big)$ occurs in the transformation behaviour of the logarithm of the Dedekind eta-function under substitutions from the modular group. In 1892, he showed the following reciprocity relation for Dedekind sums:
\begin{equation*}
S(h,m)+S(m,h)=\frac{1}{12}\bigg(\frac{h}{m}+\frac{1}{hm}+\frac{m}{h}\bigg)-\frac{1}{4},
\end{equation*}
if $h$ and $m$ are relatively prime positive integers. \par
Apostol considered the generalized Dedekind sums which are given by
\begin{equation}\label{12-1}
S_p(h,m)= \sum_{\mu=1}^{m-1}\frac{\mu}{m}\overline{B}_{p}\big(\frac{h\mu}{m}\big),
\end{equation}
and showed in [1] that they satisfy the reciprocity relation
\begin{align*}
(p+1)&\big(hm^pS_p(h,m) + mh^pS_p(m,h)\big) \\
&=pB_{p+1}+\sum_{s=0}^{p+1}\binom{p+1}{s}(-1)^sB_sB_{p+1-s}h^sm^{p+1-s}.
\end{align*} \par
In this paper, we consider the poly-Dedekind sums which are defined by
\begin{equation*}
S_p^{(k)}(h,m)= \sum_{\mu=1}^{m-1}\frac{\mu}{m}\overline{B}_{p}^{(k)}\bigg(\frac{h\mu}{m}\bigg),
\end{equation*}
where $B_{p}^{(k)}(x)$ are the type 2 poly-Bernoulli polynomials of index $k$ (see \eqref{9}) and $\overline{B}_p^{(k)}(x)=B_{p}^{(k)}(\langle x \rangle)$ are the type 2 poly-Bernoulli functions of index $k$. Note here that $S_p^{(1)}(h,m)=S_p(h,m)$. We show the following reciprocity relation for the poly-Dedekind sums given by (see Theorem 10)
\begin{align*}
&hm^{p}S_{p}^{(k)}(h,m)+mh^{p}S_{p}^{(k)}(m,h)\\
&=\sum_{\mu=0}^{m-1}\sum_{j=0}^{p}\sum_{\nu=0}^{h-1}\sum_{l=1}^{p-j+1}\frac{(mh)^{j-1}\binom{p}{j}S_{1}(p-j+1,l)}{(p-j+1)l^{k-1}}\big((\mu h)m^{p-j}+(m\nu) h^{p-j}\big)\overline{B}_{j}\bigg(\frac{\nu}{h}+\frac{\mu}{m}\bigg).
\end{align*}
For $k=1$, this reciprocity relation for the poly-Dedekind sums reduces to that for the generalized Dedekind sums given by (see Corollary 11)
\begin{align*}
&hm^{p}S_{p}(h,m)+mh^{p}S_{p}(m,h)\\
&=\ \sum_{\nu=0}^{m-1}\sum_{\nu=0}^{h-1}(mh)^{p-1}(\mu h+m\nu)\overline{B}_{p}\bigg(\frac{\nu}{h}+\frac{\mu}{m}\bigg).\\
\end{align*} \par
In Section 2, we will derive various facts about the type 2 poly-Bernoulli polynomials that will be needed in the next section. In Section 3, we will define the poly-Dedekind sums and demonstrate a reciprocity relation for them.

\section{\bf On type 2 poly-Bernoulli polynomials}
From \eqref{9}, we note that
\begin{align}
\frac{\mathrm{Ei}_{k}(\log(1+t))}{e^{t}-1}e^{xt}\ &=\ \sum_{l=0}^{\infty}B_{l}^{(k)}\frac{t^{l}}{l!}\sum_{m=0}^{\infty}\frac{x^{m}}{m!}t^{m} \label{13} \\
&=\ \sum_{n=0}^{\infty}\bigg(\sum_{l=0}^{n}\binom{n}{l}B_{l}^{(k)}x^{n-l}\bigg)\frac{t^{n}}{n!}. \nonumber
\end{align}
Thus, by \eqref{13}, we get
\begin{align}
B_{n}^{(k)}(x)\ &=\ \sum_{l=0}^{n}\binom{n}{l}B_{l}^{(k)}x^{n-l},\quad (n\ge 0).\label{14}
\end{align}
\par
By \eqref{14}, we get
\begin{equation}
\frac{d}{dx}B_{n}^{(k)}(x)=nB_{n-1}^{(k)}(x),\quad(n\ge 1).\label{15}
\end{equation}
From \eqref{9}, we have
\begin{align}
\mathrm{Ei}_{k}\big(\log(1+t)\big)\ &=\ \sum_{l=0}^{\infty}B_{l}^{(k)}\frac{t^{l}}{l!}\big(e^{t}-1\big) \label{16} \\
&=\ \sum_{n=0}^{\infty}\big(B_{n}^{(k)}(1)-B_{n}^{(k)}\big)\frac{t^{n}}{n!}\ =\ \sum_{n=1}^{\infty}\big(B_{n}^{(k)}(1)-B_{n}^{(k)}\big)\frac{t^{n}}{n!}.  \nonumber
\end{align}
On the other hand,
\begin{align}
\mathrm{Ei}_{k}\big(\log(1+t)\big)\ &=\ \sum_{m=1}^{\infty}\frac{\big(\log(1+t)\big)^{m}}{m^{k}(m-1)!}\ =\ \sum_{m=1}^{\infty}\frac{1}{m^{k-1}}\frac{1}{m!}\big(\log(1+t)\big)^{m} \label{17} \\
&=\ \sum_{m=1}^{\infty}\frac{1}{m^{k-1}}\sum_{n=m}^{\infty}S_{1}(n,m)\frac{t^{n}}{n!} \nonumber \\
&=\ \sum_{n=1}^{\infty}\bigg(\sum_{m=1}^{n}\frac{1}{m^{k-1}}S_{1}(n,m)\bigg)\frac{t^{n}}{n!},\nonumber
\end{align}
where $S_{1}(n,m)$ are the Stirling numbers of the first kind. \par
Therefore, by \eqref{16} and \eqref{17}, we obtain the following theorem.
\begin{theorem}
For $n\ge 1$, we have
\begin{displaymath}
B_{n}^{(k)}(1)-B_{n}^{(k)}=\sum_{m=1}^{n}S_{1}(n,m)\frac{1}{m^{k-1}},\quad (k\in\mathbb{Z}).
\end{displaymath}
\end{theorem}
By Theorem 2, we get
\begin{displaymath}
B_{n}^{(1)}-B_{n}^{(1)}=\delta_{1,n}, \ B_{0}^{(k)}=1,\ B_{1}^{(k)}=-1+\frac{1}{2^{k}},\dots,
\end{displaymath}
where $\delta_{n,k}$ is the Kronecker's symbol. \par

With \eqref{15} in mind, we now compute
\begin{align}
\bigg(\frac{d}{dx}\bigg)^{s-1}\big(xB_{p}^{(k)}(x)\big)\bigg|_{x=1}\ &=\ \sum_{l=0}^{s-1}\binom{s-1}{l}\bigg(\bigg(\frac{d}{dx}\bigg)^{l}x\bigg)\bigg(\bigg(\frac{d}{dx}\bigg)^{s-1-l}B_{p}^{(k)}(x)\bigg)\bigg|_{x=1} \label{18} \\
&=\ \bigg(\frac{d}{dx}\bigg)^{s-1}B_{p}^{(k)}(x)\bigg|_{x=1}+\binom{s-1}{1}\bigg(\frac{d}{dx}\bigg)^{s-2}B_{p}^{(k)}(x)\bigg|_{x=1} \nonumber \\
&=\ \frac{s!}{p+1}\binom{p+1}{s}B_{p-s+1}^{(k)}(1) + \frac{(s-1)s!}{(p+1)(p+2)}\binom{p+2}{s}B_{p-s+2}^{(k)}(1).\nonumber
\end{align}
On the other hand, by \eqref{14}, we get
\begin{align}
\bigg(\frac{d}{dx}\bigg)^{s-1}\big(xB_{p}^{(k)}(x)\big)\big|_{x=1}\ &=\ \sum_{\nu=0}^{p}\binom{p}{\nu}B_{\nu}^{(k)}\bigg(\bigg(\frac{d}{dx}\bigg)^{s-1}x^{p-\nu+1}\bigg)\bigg|_{x=1} \label{19} \\
&=\ \sum_{\nu=0}^{p}\binom{p}{\nu}B_{\nu}^{(k)}(p-\nu+1)\cdots(p-\nu-s+3) \nonumber \\
&=\ \sum_{\nu=0}^{p}\binom{p}{\nu}\frac{s!B_{\nu}^{(k)}}{p-\nu+2}\binom{p-\nu+2}{s}. \nonumber
\end{align}
Therefore, by \eqref{18} and \eqref{19}, we obtain the following theorem.
\begin{theorem}
For $s,p\in\mathbb{N}$, we have
\begin{displaymath}
\sum_{\nu=0}^{p}\binom{p}{\nu}\binom{p-\nu+2}{s}\frac{B_{\nu}^{(k)}}{p-\nu+2}=\binom{p+1}{s}\frac{B_{p-s+1}^{(k)}(1)}{p+1}+\frac{s-1}{p+1}\binom{p+2}{s}\frac{B_{p-s+2}^{(k)}(1)}{p+2}.
\end{displaymath}
\end{theorem}
Now, we observe that
\begin{align}
&\sum_{\nu=0}^{p}\binom{p}{\nu}\binom{p-\nu+2}{s}\frac{B_{\nu}^{(k)}}{p-\nu+2}=\sum_{\nu=0}^{p-s+2}\frac{\binom{p}{\nu}\binom{p-\nu+2}{s}}{p-\nu+2} B_{\nu}^{(k)} \label{20} \\
&=\sum_{\nu=0}^{p-s+1}\frac{\binom{p}{\nu}\binom{p-\nu+2}{s}}{p-\nu+2}B_{\nu}^{(k)}+\frac{1}{s}\binom{p}{s-2}B_{p-s+2}^{(k)}. \nonumber
\end{align}
Therefore, by Theorem 3 and \eqref{20}, we obtain the following corollary.
\begin{corollary}
For $s,p\in\mathbb{N}$, we have
\begin{align*}
\sum_{\nu=0}^{p-s+1}&\binom{p}{\nu}\binom{p-\nu+2}{s}\frac{B_{\nu}^{(k)}}{p-\nu+2}\\
&=\binom{p+1}{s}\frac{B_{p-s+1}^{(k)}(1)}{p+1}+\frac{s-1}{p+1}\binom{p+2}{s}\frac{B_{p-s+2}^{(k)}(1)}{p+2}-\frac{1}{s}\binom{p}{s-2}B_{p-s+2}^{(k)}.
\end{align*}
\end{corollary}
From \eqref{15}, we have
\begin{align}
\int_{0}^{1}xB_{p}^{(k)}(x)dx\ &=\ \bigg[x\frac{B_{p+1}^{(k)}(x)}{p+1}\bigg]_{0}^{1}-\frac{1}{p+1}\int_{0}^{1}B_{p+1}^{(k)}(x)dx \label{21} \\
&=\ \frac{B_{p+1}^{(k)}(1)}{p+1}-\frac{1}{p+1}\bigg[\frac{1}{p+2}B_{p+2}^{(k)}(x)\bigg]_{0}^{1}\nonumber \\
&=\ \frac{B_{p+1}^{(k)}(1)}{p+1}-\frac{B_{p+2}^{(k)}(1)}{(p+1)(p+2)}+\frac{B_{p+2}^{(k)}}{(p+1)(p+2)}.\nonumber
\end{align}
On the other hand, by \eqref{14}, we get
\begin{align}
\int_{0}^{1}xB_{p}^{(k)}(x)dx\ &=\ \sum_{s=0}^{p}\binom{p}{s}B_{s}^{(k)}\int_{0}^{1}x^{p-s+1}dx\label{22}\\
&=\ \sum_{s=0}^{p}\binom{p}{s}B_{s}^{(k)}\frac{1}{p+2-s}.\nonumber
\end{align}
Therefore, by \eqref{21} and \eqref{22}, we obtain the following theorem.
\begin{theorem}
For $p\in\mathbb{N}$, we have
\begin{displaymath}
\sum_{s=0}^{p}\binom{p}{s}B_{s}^{(k)}\frac{1}{p+2-s}=\frac{B_{p+1}^{(k)}(1)}{p+1}-\frac{B_{p+2}^{(k)}(1)}{(p+1)(p+2)}+\frac{B_{p+2}^{(k)}}{(p+1)(p+2)}.
\end{displaymath}
\end{theorem}

\section{\bf Poly-Dedekind Sums}
Apostol considered the generalized Dedekind sums which are given by
\begin{equation}
S_{p}(h,m)=\sum_{\mu=1}^{m-1}(\mu/m)\overline{B}_{p}(h\mu/m),\quad(h,m,p\in\mathbb{N}), \label{23}
\end{equation}
where $\overline{B}_{p}(h\mu/m)=B_{p}\big(\langle h\mu/m\rangle\big)$. \par
Note that, for any relatively prime positive integers $h,m$, we have
\begin{align*}
S_{1}(h,m)\ &=\ \sum_{\mu=1}^{m-1}(\mu/m)\overline{B}_{1}(h\mu/m) \\
&=\ \sum_{\mu=1}^{m-1}(\!(\mu/m)\!)(\!(h\mu/m)\!)\ =\ S(h,m).
\end{align*}
In this section, we consider the {\it{poly-Dedekind sums}} which are given by
\begin{equation}
S_{p}^{(k)}(h,m)\ =\ \sum_{\mu=1}^{m-1}\big(\mu/m\big)\overline{B}_{p}^{(k)}\big(h\mu/m\big),\label{24}
\end{equation}
where $h,m,p\in\mathbb{N}, k\in\mathbb{Z}$, and $\overline{B}_p^{(k)}(x)=B_p^{(k)}(\langle x \rangle)$ are the type 2 poly-Bernoulli functions of index $k$. \par
Note that
\begin{displaymath}
S_{p}^{(1)}(h,m)\ =\ \sum_{\mu=1}^{m-1}\big(\mu/m\big)\overline{B}_{p}\big(h\mu/m\big)=S_{p}(h,m).
\end{displaymath}
Assume now that $h=1$. Then we have
\begin{align}
S_{p}^{(k)}(1,m)\ &=\ \sum_{\mu=1}^{m-1}\big(\mu/m\big)\overline{B}_{p}^{(k)}\big(\mu/m\big) \label{25}\\
    &=\ \sum_{\mu=1}^{m-1}\big(\mu/m)\sum_{\nu=0}^{p}\binom{p}{\nu}B_{\nu}^{(k)}\big(\mu/m\big)^{p-\nu}\nonumber \\
    &=\ \sum_{\nu=0}^{p}\binom{p}{\nu}B_{\nu}^{(k)}m^{-(p-\nu+1)}\sum_{\mu=1}^{m-1}\mu^{p+1-\nu}\nonumber\\
    &=\ \sum_{\nu=0}^{p}\binom{p}{\nu}B_{\nu}^{(k)}m^{-(p+1-\nu)}\frac{1}{p+2-\nu}\big(B_{p+2-\nu}(m)-B_{p+2-\nu}\big). \nonumber
\end{align}
From \eqref{5}, we have
\begin{align}
B_{p+2-\nu}(m)-B_{p+2-\nu}\ &=\ \sum_{i=0}^{p+2-\nu}\binom{p+2-\nu}{i}B_{i}m^{p+2-\nu-i}-B_{p+2-\nu} \label{26} \\
&=\ \sum_{i=0}^{p+1-\nu} \binom{p+2-\nu}{i} B_{i}m^{p+2-\nu-i}. \nonumber
\end{align}
By \eqref{25} and \eqref{26}, we get
\begin{align}
S_{p}^{(k)}(1,m)\ &=\ \sum_{\nu=0}^{p}\binom{p}{\nu}B_{\nu}^{(k)}m^{-(p+1-\nu)}\frac{1}{p+2-\nu}\sum_{i=0}^{p+1-\nu}\binom{p+2-\nu}{i}B_{i}m^{p+2-\nu-i}\label{27} \\
&=\     \ \frac{1}{m^{p}}\sum_{\nu=0}^{p}\binom{p}{\nu}\frac{B_{\nu}^{(k)}}{p+2-\nu}\sum_{i=0}^{p+1-\nu}\binom{p+2-\nu}{i}B_{i}m^{p+1-i}.\nonumber
\end{align}

Now, we assume that $p$ is an odd positive integer $\ge 3$, so that $B_p=0$.
Then we have
\begin{align}
&m^{p}S_{p}^{(k)}(1,m)\ =\ \sum_{\nu=0}^{p}\binom{p}{\nu}\frac{B_{\nu}^{(k)}}{p+2-\nu}\sum_{i=0}^{p+1-\nu}\binom{p+2-\nu}{i}B_{i}m^{p+1-i}\label{28}\\
&=\ \sum_{\nu=0}^{p}\binom{p}{\nu}\frac{B_{\nu}^{(k)}}{p+2-\nu}m^{p+1}+\sum_{\nu=0}^{p}\binom{p}{\nu}\frac{B_{\nu}^{(k)}}{p+2-\nu}\sum_{i=1}^{p+1-\nu}\binom{p+2-\nu}{i}B_{i}m^{p+1-i}\nonumber     \\
&=\ \sum_{\nu=0}^{p}\binom{p}{\nu}\frac{B_{\nu}^{(k)}}{p+2-\nu}m^{p+1}+\sum_{i=1}^{p+1}\sum_{\nu=0}^{p+1-i}\binom{p}{\nu} \binom{p+2-\nu}{i}\frac{B_{\nu}^{(k)}}{p+2-\nu}B_{i}m^{p+1-i}\nonumber \\
&=\ \sum_{\nu=0}^{p}\binom{p}{\nu}\frac{B_{\nu}^{(k)}}{p+2-\nu}m^{p+1}+\sum_{i=1}^{p-1}\sum_{\nu=0}^{p+1-i}\binom{p}{\nu} \binom{p+2-\nu}{i}\frac{B_{\nu}^{(k)}}{p+2-\nu} B_{i} m^{p+1-i}\nonumber  \\
&\quad +\frac{1}{p+2}\binom{p+2}{p+1}B_{p+1}+\sum_{\nu=0}^{1}\binom{p}{\nu}\binom{p+2-\nu}{p}\frac{B_{\nu}^{(k)}}{p+2-\nu}B_{p}m \nonumber \\
&=\ \sum_{\nu=0}^{p}\binom{p}{\nu}\frac{B_{\nu}^{(k)}}{p+2-\nu}m^{p+1}+ \sum_{i=1}^{p-1}\sum_{\nu=0}^{p+1-i}\binom{p}{\nu}\frac{\binom{p+2-\nu}{i} }{p+2-\nu}B_{\nu}^{(k)}B_{i}m^{p+1-i}+B_{p+1}\nonumber.
\end{align}
Therefore, by \eqref{28}, we obtain the following proposition.
\begin{proposition}
Let $p$ be an odd positive integer $\ge 3$. Then we have
\begin{displaymath}
m^{p}S_{p}^{(k)}(1,m)= \sum_{\nu=0}^{p}\binom{p}{\nu}\frac{B_{\nu}^{(k)}}{p+2-\nu}m^{p+1}+ \sum_{i=1}^{p-1}\sum_{\nu=0}^{p+1-i}\binom{p}{\nu} \binom{p+2-\nu}{i}\frac{B_{\nu}^{(k)}}{p+2-\nu}B_{i}m^{p+1-i}+B_{p+1}.
\end{displaymath}
\end{proposition}
We still assume that $p$ is an odd positive integer $\ge 3$, so that $B_p=0$.
Then, from Corollary 4, Theorem 5 and Proposition 6, we note that
\begin{align}
&m^{p}S_{p}^{(k)}(1,m)\label{29} \\
&=\ \sum_{\nu=0}^{p}\binom{p}{\nu}\frac{B_{\nu}^{(k)}}{p+2-\nu}m^{p+1}+ \sum_{i=1}^{p-1}\sum_{\nu=0}^{p+1-i}\binom{p}{\nu} \binom{p+2-\nu}{i}\frac{B_{\nu}^{(k)}}{p+2-\nu}B_{i}m^{p+1-i}+B_{p+1}.\nonumber \\
&=\bigg(\frac{B_{p+1}^{(k)}(1)}{p+1}-\frac{B_{p+2}^{(k)}(1)}{(p+1)(p+2)}+\frac{B_{p+2}^{(k)}}{(p+1)(p+2)}\bigg)m^{p+1}+B_{p+1}\nonumber \\
&+\sum_{i=1}^{p-1}\bigg(\binom{p+1}{i}\frac{B_{p+1-i}^{(k)}(1)}{p+1}+\frac{(i-1)}{(p+1)(p+2)}\binom{p+2}{i}B_{p+2-i}^{(k)}(1)-\binom{p}{i-2}\frac{1}{i}B_{p+2-i}^{(k)}\bigg)B_{i}m^{p+1-i}.\nonumber
\end{align}
To proceed further, we note that $\binom{p}{i-2}\frac{p+1}{i}=\frac{1}{p+2}\binom{p+2}{i}(i-1)$, for $i \ge 1$, and that $B_1^{(k)}(1)-B_1^{(k)}=1$, from Theorem 2.
Then, from \eqref{29}, we see that
\begin{align}
&(p+1)m^{p}S_{p}^{(k)}(1,m)=\bigg(B_{p+1}^{(k)}(1)-\frac{B_{p+2}^{(k)}(1)}{p+2}+\frac{B_{p+2}^{(k)}}{p+2}\bigg)m^{p+1} \label{30} \\
&\quad+\sum_{i=1}^{p-1}\binom{p+1}{i}B_{i}B_{p+1-i}^{(k)}(1)m^{p+1-i}+(p+1)B_{p+1} \nonumber \\
&\quad+\frac{1}{p+2}\sum_{i=1}^{p-1}\binom{p+2}{i}(i-1)B_{i}B_{p+2-i}^{(k)}(1)m^{p+1-i}-\sum_{i=1}^{p-1}\binom{p}{i-2}\frac{(p+1)}{i}B^{(k)}_{p+2-i}B_{i}m^{p+1-i}
\nonumber \\
&=m^{p+1}B_{p+1}^{(k)}(1)+\sum_{i=1}^{p-1}\binom{p+1}{i}B_{i}m^{p+1-i}B_{p+1-i}^{(k)}(1)+B_{p+1}\nonumber \\
&+\frac{1}{p+2}(-1)m^{p+1}\big(B_{p+2}^{(k)}(1)-B_{p+2}^{(k)}\big)+\frac{1}{p+2}\sum_{i=1}^{p-1}\binom{p+2}{i}(i-1)B_{i}m^{p+1-i} \big(B_{p+2-i}^{(k)}(1)-B_{p+2-i}^{(k)}\big)
\nonumber \\
&+pB_{p+1} \nonumber \\
&= \sum_{i=0}^{p+1}\binom{p+1}{i}B_{i}m^{p+1-i}B_{p+1-i}^{(k)}(1)+\frac{1}{p+2}\sum_{i=0}^{p+1}\binom{p+2}{i}(i-1)B_{i}m^{p+1-i} \big(B_{p+2-i}^{(k)}(1)-B_{p+2-i}^{(k)}\big). \nonumber
\end{align}
Therefore, by \eqref{30}, we obtain the following theorem. \begin{theorem}
For $m\in\mathbb{N}$, and any odd positive integer $ p \ge 3$ , we have
\begin{align*}
&(p+1)m^{p}S_{p}^{(k)}(1,m) \\
&= \sum_{i=0}^{p+1}\binom{p+1}{i}B_{i}m^{p+1-i}B_{p+1-i}^{(k)}(1)+\frac{1}{p+2}\sum_{i=0}^{p+1}\binom{p+2}{i}(i-1)B_{i}m^{p+1-i} \big(B_{p+2-i}^{(k)}(1)-B_{p+2-i}^{(k)}\big).
\end{align*}
\end{theorem}
Now we employ the symbolic notation as
\begin{displaymath}
B_{n}(x)=(B+x)^{n},\quad B_{n}^{(k)}(x)=\big(B^{(k)}+x
\big)^{n},\quad (n\ge 0).
\end{displaymath}
Assume that $h, m$ are relatively prime positive integers. Then we see that
\begin{align}
    &m^{p}\sum_{\mu=0}^{m-1}\sum_{s=0}^{p+1}\binom{p+1}{s}h^{s}B_{s}^{(k)}\big(\mu/m\big)B_{p+1-s}\big(h-[h\mu/m]\big)\label{31} \\
    &=m^{p} \sum_{\mu=0}^{m-1}\sum_{s=0}^{p+1}\binom{p+1}{s}h^{s}\big(B^{(k)}+\mu m^{-1}\big)^{s}\big(B+h-[h\mu/m]\big)^{p+1-s} \nonumber\\
    &=m^{p}\sum_{\mu=0}^{m-1}\bigg(hB^{(k)}+h\mu m^{-1}+B+h-[h\mu/m]\bigg)^{p+1} \nonumber\\
    &=m^{p}\sum_{\mu=0}^{m-1}\bigg(hB^{(k)}+h+B+\frac{1}{2}+\frac{h\mu}{m} -[h\mu/m]-\frac{1}{2}\bigg)^{p+1}\nonumber \\
    &=m^{p}\sum_{\mu=0}^{m-1}\bigg(hB^{(k)}+h+B+\frac{1}{2}+\overline{B}_{1}\big(h\mu/m\big)\bigg)^{p+1}.\nonumber
\end{align}
Now, as the index $\mu$ ranges over the values $\mu=0,1,2,\dots,m-1$, the product $h\mu$ ranges over a complete residue system modulo $m$ and due to the periodicity of $\overline{B}_{1}(x)$, the term $\overline{B}_{1}(h\mu,m)$ may be replaced by $\overline{B}_{1}(\mu/m)$, without altering the sum over $\mu$. Thus the sum \eqref{31} is equal to
\begin{align}
&\ m^{p}\sum_{m=0}^{m-1}\bigg(hB^{(k)}+h+B+\frac{1}{2}+\overline{B}_{1}\bigg(\frac{\mu}{m}\bigg)\bigg)^{p+1}\label{32}\\
&=\ m^{p}\sum_{m=0}^{m-1}\bigg(h\big(B^{(k)}+1\big)+B+\frac{\mu}{m}\bigg)\bigg)^{p+1}\nonumber \\
&=\ m^{p}\sum_{\mu=0}^{m-1}\sum_{s=0}^{p+1}\binom{p+1}{s}\bigg(B+\frac{\mu}{m}\bigg)^{s}h^{p+1-s}\big(B^{(k)}+1\big)^{p+1-s}\nonumber \\
&=\ m^{p}\sum_{\mu=0}^{m-1}\sum_{s=0}^{p+1}\binom{p+1}{s}B_{s}\bigg(\frac{\mu}{m}\bigg)h^{p+1-s}B^{(k)}_{p+1-s}(1)\nonumber \\
&=\ \sum_{s=0}^{p+1}\binom{p+1}{s}m^{s-1}\sum_{\mu=0}^{m-1}B_{s}\bigg(\frac{\mu}{m}\bigg)(mh)^{p+1-s}B_{p+1-s}^{(k)}(1)\nonumber \\
&=\ \sum_{s=0}^{p+1}\binom{p+1}{s}B_{s}(mh)^{p+1-s}B_{p+1-s}^{(k)}(1),\nonumber
\end{align}
where we used the fact (a) in Lemma 1.\par
Therefore, we obtain the following theorem.
\begin{theorem}
For $m,n,h\in\mathbb{N}$, with $(h,m)=1$, and $p$ is any positive odd integer $ \ge 3$, we have
\begin{align*}
\sum_{s=0}^{p+1}\binom{p+1}{s}B_{s}B_{p+1-s}^{(k)}(1)(mh)^{p+1-s}=m^{p}\sum_{\mu=0}^{m-1}\sum_{s=0}^{p+1}\binom{p+1}{s}h^{s}B_{s}^{(k)}\big(\mu/m\big)B_{p+1-s}\bigg(h-\bigg[\frac{h\mu}{m}\bigg]\bigg).
\end{align*}
\end{theorem}
Now, we observe that
\begin{align}
\sum_{n=0}^{\infty}B_{n}^{(k)}(x)\frac{t^{n}}{n!}\ &=\ \frac{\mathrm{Ei}_{k}\big(\log(1+t)\big)}{e^{t}-1}e^{xt}=\frac{\mathrm{Ei}_{k}\big(\log(1+t)\big)}{e^{dt}-1}\sum_{i=0}^{d-1}e^{(i+x)t}\label{33}\\
&=\ \frac{\mathrm{Ei}_{k}\big(\log(1+t)\big)}{dt}\sum_{i=0}^{d-1}e^{(i+x)t}\frac{dt}{e^{dt}-1}\nonumber\\
&=\ \sum_{j=0}^{\infty}d^{j-1}\sum_{i=0}^{d-1}B_{j}\bigg(\frac{x+i}{d}\bigg)\frac{t^{j}}{j!}\frac{1}{t}\sum_{l=1}^{\infty}\frac{\big(\log(1+t)\big)^{l}}{(l-1)!l^{k}}\nonumber \\
&=\ \sum_{j=0}^{\infty}d^{j-1}\sum_{i=0}^{d-1}B_{j}\bigg(\frac{x+i}{d}\bigg)\frac{t^{j}}{j!}\frac{1}{t}\sum_{l=1}^{\infty}\frac{1}{l^{k-1}}\sum_{m=l}^{\infty}S_{1}(m,l)\frac{t^{m}}{m!}\nonumber\\
&=\ \sum_{j=0}^{\infty}d^{j-1}\sum_{i=0}^{d-1}B_{j}\bigg(\frac{x+i}{d}\bigg)\frac{t^{j}}{j!}\frac{1}{t}\sum_{m=1}^{\infty}\sum_{l=1}^{m}\frac{S_{1}(m,l)}{l^{k-1}}\frac{t^{m}}{m!}\nonumber\\
&=\ \sum_{j=0}^{\infty}d^{j-1}\sum_{i=0}^{d-1}B_{j}\bigg(\frac{x+i}{d}\bigg)\frac{t^{j}}{j!}\sum_{m=0}^{\infty}\sum_{l=1}^{m+1}\frac{S_{1}(m+1,l)}{l^{k-1}(m+1)}\frac{t^{m}}{m!}\nonumber\\
&=\ \sum_{n=0}^{\infty}\bigg(\sum_{j=0}^{n}\sum_{i=0}^{d-1}\sum_{l=1}^{n-j+1}\binom{n}{j}d^{j-1}B_{j}\bigg(\frac{x+i}{d}\bigg)\frac{S_{1}(n-j+1,l)}{(n-j+1)l^{k-1}}\bigg)\frac{t^{n}}{n!}, \nonumber
\end{align}
where $d$ is a positive integer. \par
Therefore, by comparing the coefficients on both sides of \eqref{33}, we obtain the following theorem.
\begin{theorem}
For $k\in\mathbb{Z}$, $d\in\mathbb{N}$ and $n\ge 0$, we have
\begin{displaymath}
B_{n}^{(k)}(x)=\sum_{j=0}^{n}\sum_{i=0}^{d-1}\sum_{l=1}^{n-j+1}\binom{n}{j}d^{j-1}B_{j}\bigg(\frac{x+i}{d}\bigg)\frac{S_{1}(n-j+1,l)}{(n-j+1)l^{k-1}}.
\end{displaymath}
\end{theorem}

From \eqref{24}, and by using Theorem 9 and (c) in Lemma 1, we see that
\begin{align}
& hm^{p}S_{p}^{(k)}(h,m)+mh^{p}S_{p}^{(k)}(m,h)\label{34}\\
&=hm^{p}\sum_{\mu=0}^{m-1}\frac{\mu}{m}\overline{B}_{p}^{(k)}\bigg(\frac{h\mu}{m}\bigg)+mh^{p}\sum_{\nu=0}^{h-1}\frac{\nu}{h}\overline{B}_{p}^{(k)}\bigg(\frac{m\nu}{h}\bigg) \nonumber\\
&=hm^{p}\sum_{\mu=0}^{m-1}\frac{\mu}{m}\sum_{j=0}^{p}h^{j-1}\binom{p}{j}\sum_{\nu=0}^{h-1}\sum_{l=1}^{p-j+1}\frac{S_{1}(p-j+1,l)}{(p-j+1)l^{k-1}}\overline{B}_{j}\bigg(\frac{\mu}{m}+\frac{\nu}{h}\bigg)\nonumber \\
&+mh^{p}\sum_{\nu=0}^{h-1}\frac{\nu}{h}\sum_{j=0}^{p}m^{j-1}\binom{p}{j}\sum_{\mu=0}^{m-1}\sum_{l=1}^{p-j+1} \frac{S_{1}(p-j+1,l)}{(p-j+1)l^{k-1}}\overline{B}_{j}\bigg(\frac{\nu}{h}+\frac{\mu}{m}\bigg)\nonumber \\
&=\sum_{\mu=0}^{m-1}\frac{\mu}{m}\sum_{j=0}^{p}m^{p-j}(mh)^{j}\binom{p}{j}\sum_{\nu=0}^{h-1}\sum_{l=1}^{p-j+1}\overline{B}_{j}\bigg(\frac{\mu}{m}+\frac{\nu}{h}\bigg)\frac{S_{1}(p-j+1,l)}{(p-j+1)l^{k-1}}\nonumber \\
&+\sum_{\nu=0}^{h-1}\frac{\nu}{h}\sum_{j=0}^{p}h^{p-j}(mh)^{j}\binom{p}{j}\sum_{\mu=0}^{m-1}\sum_{l=1}^{p-j+1}\overline{B}_{j}\bigg(\frac{\nu}{h}+\frac{\mu}{m}\bigg)\frac{S_{1}(p-j+1,l)}{(p-j+1)l^{k-1}}\nonumber \\
&=\sum_{\mu=0}^{m-1}\sum_{j=0}^{p}\sum_{\nu=0}^{h-1}\sum_{l=1}^{p-j+1}(\mu h)(mh)^{-1}m^{p-j}(mh)^{j}\binom{p}{j}\overline{B}_{j}\bigg(\frac{\mu}{m}+\frac{\nu}{h}\bigg)\frac{S_{1}(p-j+1,l)}{(p-j+1)l^{k-1}}\nonumber \\
&+\sum_{\mu=0}^{m-1}\sum_{j=0}^{p}\sum_{\nu=0}^{h-1}\sum_{l=1}^{p-j+1}(m\nu)(mh)^{-1}h^{p-j}(mh)^{j}\binom{p}{j}\overline{B}_{j}\bigg(\frac{\nu}{h}+\frac{\mu}{m}\bigg)\frac{S_{1}(p-j+1,l)}{(p-j+1)l^{k-1}}\nonumber \\
&=\sum_{\mu=0}^{m-1}\sum_{j=0}^{p}\sum_{\nu=0}^{h-1}\sum_{l=1}^{p-j+1}\frac{(mh)^{j-1}\binom{p}{j}S_{1}(p-j+1,l)}{(p-j+1)l^{k-1}}\big((\mu h)m^{p-j}+(m\nu) h^{p-j}\big)\overline{B}_{j}\bigg(\frac{\nu}{h}+\frac{\mu}{m}\bigg).\nonumber
\end{align}
Therefore, we obtain the following reciprocity relation.

\begin{theorem}
For $m,h,p\in\mathbb{N}$ and $k\in\mathbb{Z}$, we have
\begin{align*}
&hm^{p}S_{p}^{(k)}(h,m)+mh^{p}S_{p}^{(k)}(m,h)\\
&=\sum_{\mu=0}^{m-1}\sum_{j=0}^{p}\sum_{\nu=0}^{h-1}\sum_{l=1}^{p-j+1}\frac{(mh)^{j-1}\binom{p}{j}S_{1}(p-j+1,l)}{(p-j+1)l^{k-1}}\big((\mu h)m^{p-j}+(m\nu) h^{p-j}\big)\overline{B}_{j}\bigg(\frac{\nu}{h}+\frac{\mu}{m}\bigg).
\end{align*}
\end{theorem}

In case of $k=1$, we obtain the following reciprocity relation for the generalized Dedekind sum defined by Apostol.

\begin{corollary}
For $m,h,p\in\mathbb{N}$, we have
\begin{align*}
&hm^{p}S_{p}(h,m)+mh^{p}S_{p}(m,h)\\
&=\ \sum_{\nu=0}^{m-1}\sum_{\nu=0}^{h-1}(mh)^{p-1}(\mu h+m\nu)\overline{B}_{p}\bigg(\frac{\nu}{h}+\frac{\mu}{m}\bigg)\\
&=\ (mh)^{p} \sum_{\nu=0}^{m-1}\sum_{\nu=0}^{h-1}(mh)^{-1}(\mu h+m\nu)\overline{B}_{p}\bigg(\frac{\nu}{h}+\frac{\mu}{m}\bigg).
\end{align*}
\end{corollary}

\section{\bf Conclusion}
 Dedekind sums are defined by
 \begin{equation}
 S(h,m)=\sum_{\mu=1}^{m}\bigg(\!\bigg(\frac{\mu}{m}\bigg)\!\bigg)\bigg(\!\bigg(\frac{h\mu}{m}\bigg)\!\bigg),
 \quad(\mathrm{see}\ [1,3,6-8,11-13]).
 \end{equation}
In 1952, Apostol considered the generalized Dedekind sums and
introduced interested and important identities and theorems
related to his generalized Dedekind sums.  These Dedekind sums are
a field that has been studied by various researchers. Recently,
the modified Hardy's polyexponential function of index $k$  is
introduced by
 \begin{equation}
 \mathrm{Ei}_{k}(x)=\sum_{n=1}^{\infty}\frac{x^{n}}{n^{k}(n-1)!},\ (k\in\mathbb{Z}),\quad
 (\mathrm{see}\ [5,10]).
 \end{equation}
In [5], the type 2 poly-Bernoulli polynomials of index $k$ are
defined in terms of the polyexponential function of index $k$ as
\begin{equation}
\frac{\mathrm{Ei}_{k}(\log(1+t))}{e^{t}-1}e^{xt}=\sum_{n=0}^{\infty}B_{n}^{(k)}(x)\frac{t^{n}}{n!},\quad(k\in\mathbb{Z}).
\label{9}
\end{equation}
In this paper, we thought of the poly-Dedekind sums from the
perspective of the Apostol's generalized Dedekind sums. That is,
we considered the poly=Dedekind sums which are derived from the
type 2 poly-Bernoulli functions and polynomials.

\vspace{0.2in}

{\bf Acknowledgements}

The authors thank Jangjeon Institute for Mathematical Science for
the support of this research.

\vspace{0.1in}

{\bf Funding}

Not applicable.

\vspace{0.1in}

{\bf Availability of data and materials}

Not applicable.

\vspace{0.1in}

{\bf Ethics approval and consent to participate}

All authors reveal that there is no ethical problem in the
production of this paper.

\vspace{0.1in}

{\bf Competing interests}

The authors declare that they have no competing interests.

\vspace{0.1in}

{\bf Consent for publication}

All authors want to publish this paper in this journal.

\vspace{0.1in}

{\bf Authors contributions}

T.K. and D.S.K. conceived of the framework and structured the
whole paper; T.K. and D.S.K. wrote the paper; L.-C. J. and H.L.
checked the results of the paper; D.S.K. and T.K. completed the
revision of the article. All authors have read and agreed to the
published version of the manuscript.

\end{document}